\pdfoutput=1

\documentclass[12pt]{article}

\usepackage{graphicx}
\usepackage[a4paper]{geometry}
\usepackage[linktocpage=true]{hyperref}
\usepackage{float}
\usepackage{newfloat}

\providecommand{\keywords}[1]{
	\par \vspace{1 em} \noindent \em \parbox{\linewidth}{\textbf{keywords} ---\ #1}}

\providecommand{\msc}[1]{
	\par \vspace{-0.55 em} \noindent \parbox{\linewidth}{\textbf{2010 MSC:} \em  #1}}

\emergencystretch=\maxdimen
\begin{document}
\title{\huge The toric sections: a simple introduction}
\author{Luca Moroni - \href{https://www.lucamoroni.it}{www.lucamoroni.it}\\
\normalsize Liceo Scientifico Donatelli-Pascal\\
\normalsize Milano - Italy} \date{} \normalsize \maketitle

\setlength{\parskip}{.6 em}

\begin{abstract}
\noindent We review, from a didactic point of view, the definition of a toric section and the different shapes it can take.
We'll then discuss some properties of this curve, investigate its analogies and differences with the most renowned conic section and show how to build its general quartic equation.
\par
A curious and unexpected result was to find that, with some algebraic manipulation, a toric section can also be obtained as the intersection of a cylinder with a cone.
Finally we'll show how it is possible to construct and represent toric sections in the 3D Graphics view of Geogebra.
\par
In the article only elementary algebra is used, and the requirements to follow it are just some notion of goniometry and of tridimensional analytic geometry.
\keywords {toric sections, conic sections, Cassini's oval, Bernoulli's lemniscate, quartic, bicircular quartic, geogebra}
\msc {97G40}
\end{abstract}

\section {Overview}
The curve of intersection of a torus with a plane is called \textbf{toric section}. 
Even if both surfaces are rather simple to define and are described by rather simple equations, the toric section has a rather complicated equation and can assume rather interesting shapes.
\par
In this article we'll present some types of toric sections, including Villarceau's circles, Cassini's ovals, Bernoulli's lemniscates and the Hippopedes of Proclus.
\par
We'll then discuss some properties of these curves, and compare them with the most renowned conic sections.
\par
After that we'll derive the quartic equation that represents, in Cartesian coordinates, the toric section. 
By trying to simplify this equation we have found that, with some algebraic manipulation, a toric section can also be obtained as the intersection of a cylinder with a cone. 
\par
Finally we'll show how it is possible to represent toric sections in the 3D Graphics view of Geogebra.

\clearpage

\setlength{\parskip}{.0 em}
\footnotesize \setcounter{tocdepth}{2}
\tableofcontents
\setlength{\parskip}{.6 em}
\normalsize
\vspace*{1\baselineskip}
\begin{figure}[H]
	\centering
	\includegraphics[width=4.2in]{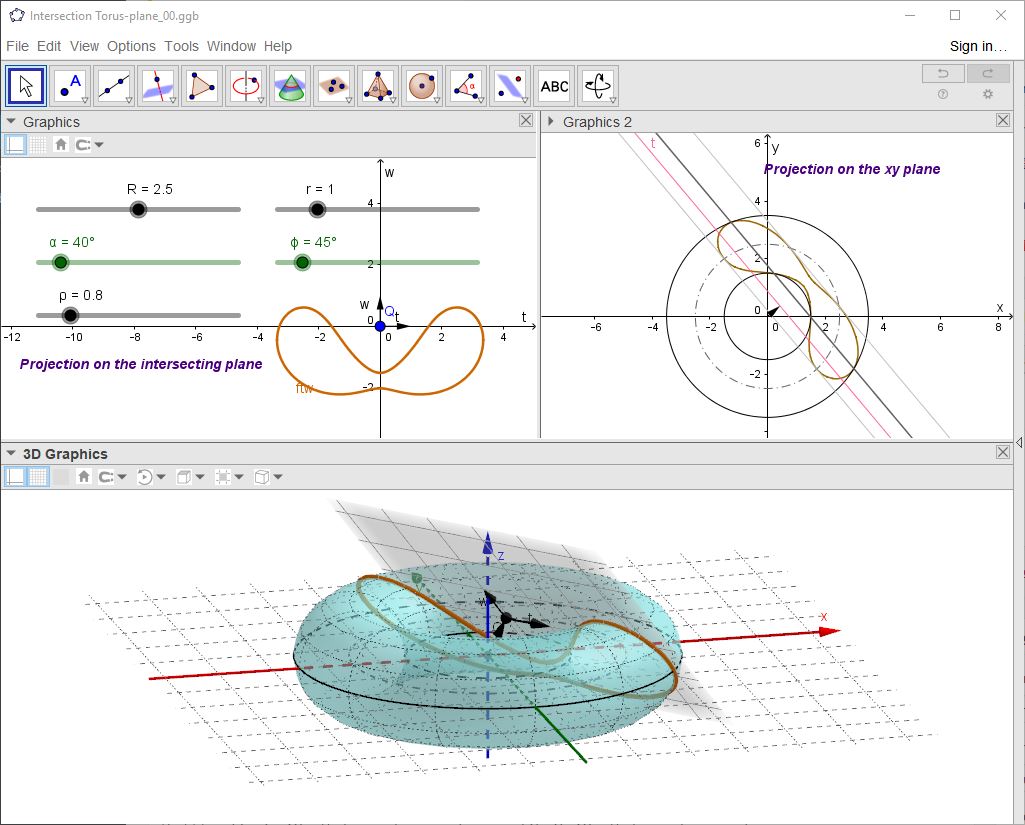}
	\caption{The torus-plane intersection simulation with Geogebra}
	\label{fig01}
\end{figure}

\clearpage

\section{The torus and the toric section}

The torus surface is generated by rotating a circle with radius $ r $ around an axis coplanar with the circle, following a second circle of radius $ R $. The two radii $ r $ and $ R $ are so the parameters that identify the torus' shape (see Figure \ref{fig10}). 

In the following lines (if not otherwise specified) we'll always assume that $R \ge r$  (to avoid self-intersections), and that, in Cartesian coordinates, the equatorial plane is the $xy$ plane and the axis of rotation is the $z$ axis.

A toric section is the analogue of a conic section as it is the intersection curve of a torus with a plane just as a conic section is the intersection curve between a conical surface and a plane. 

\renewcommand{\thefigure}{\arabic{figure}{.a}}
\begin{figure}[H]
	\centering
	\includegraphics[width=4.5in]{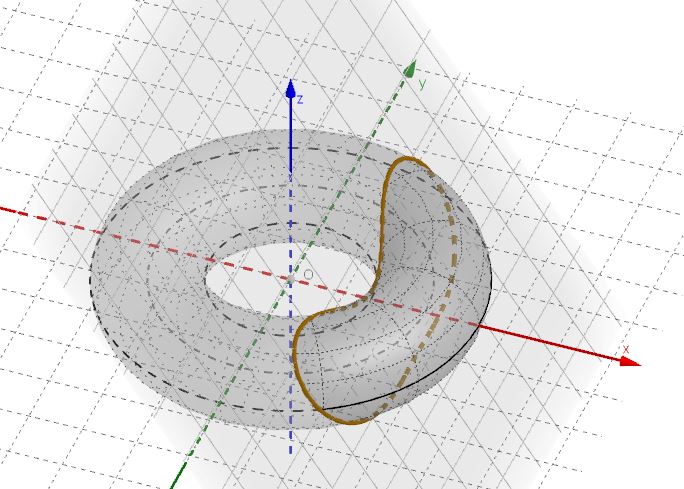}
	\caption{A toric section (red line)}
	\label{fig02a}
\end{figure}
\renewcommand{\thefigure}{\arabic{figure}{.b}}
\addtocounter {figure} {-1}
\begin{figure}[H]
	\centering
	\includegraphics[width=4.5in]{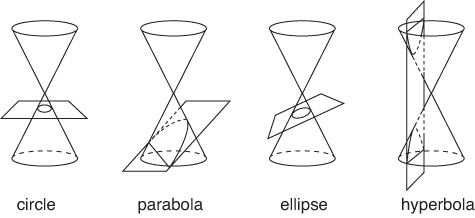}
	\caption{The \em renowned \em conic sections}
	\label{fig02b}
\end{figure}

\renewcommand{\thefigure}{\arabic{figure}}
\noindent But whilst conic sections have a (deserved) renown also due to their multiple connections with many fundamental problems of (classical) physics, toric sections are often relegated to the realm of mathematical curiosities and exotic curves. Yet...

\begin{quote} \noindent \em \footnotesize The variety of the toric sections is quite rich and the study of the case when sectioning planes are parallel to the symmetry axis of the torus dates back to antiquity. The corresponding curves are usually called "the spiric sections of Perseus" after their discoverer (circa 150 BC). Among them, one can distinguish Cassini's curves (in particular, Bernoulli's lemniscate).\em 
\par

\noindent \small Antoni Sym (2009), \emph{"Darboux's greatest love"} \cite{greatestlove}  \noindent
\end{quote}

\par

\section{The variety of toric sections}
In this section we'll shortly describe and present some particular toric sections that also have an historical relevance. Amongst them Villarceau's circles, Cassini's ovals, Bernoulli's lemniscates and the Hippopedes of Proclus.

\subsection{The central toric sections and the Villarceau's circles}
A toric section in which the cutting plane passes through the center of the torus is called \textbf{\textit{central toric section}}.
\par
\noindent Central toric sections can be circles in the following cases:
\begin{itemize}
\item the intersecting plane is the equatorial plane
\item the intersecting plane is perpendicular to the equatorial plane
\item the intersecting plane touches the torus in two isolated points\footnote{The angle formed by the intersecting plane with the equatorial plane must be $\arctan \frac{{\sqrt {{R^2} - {r^2}} }}{r}$}.
\end{itemize}
The last (less banal) case produces two circles that are called Villarceau's circles (Figure 3).

\begin{figure}[H]
	\centering
	\begin{minipage}{.49\textwidth}
	\centering
	\includegraphics[width=0.8\linewidth]{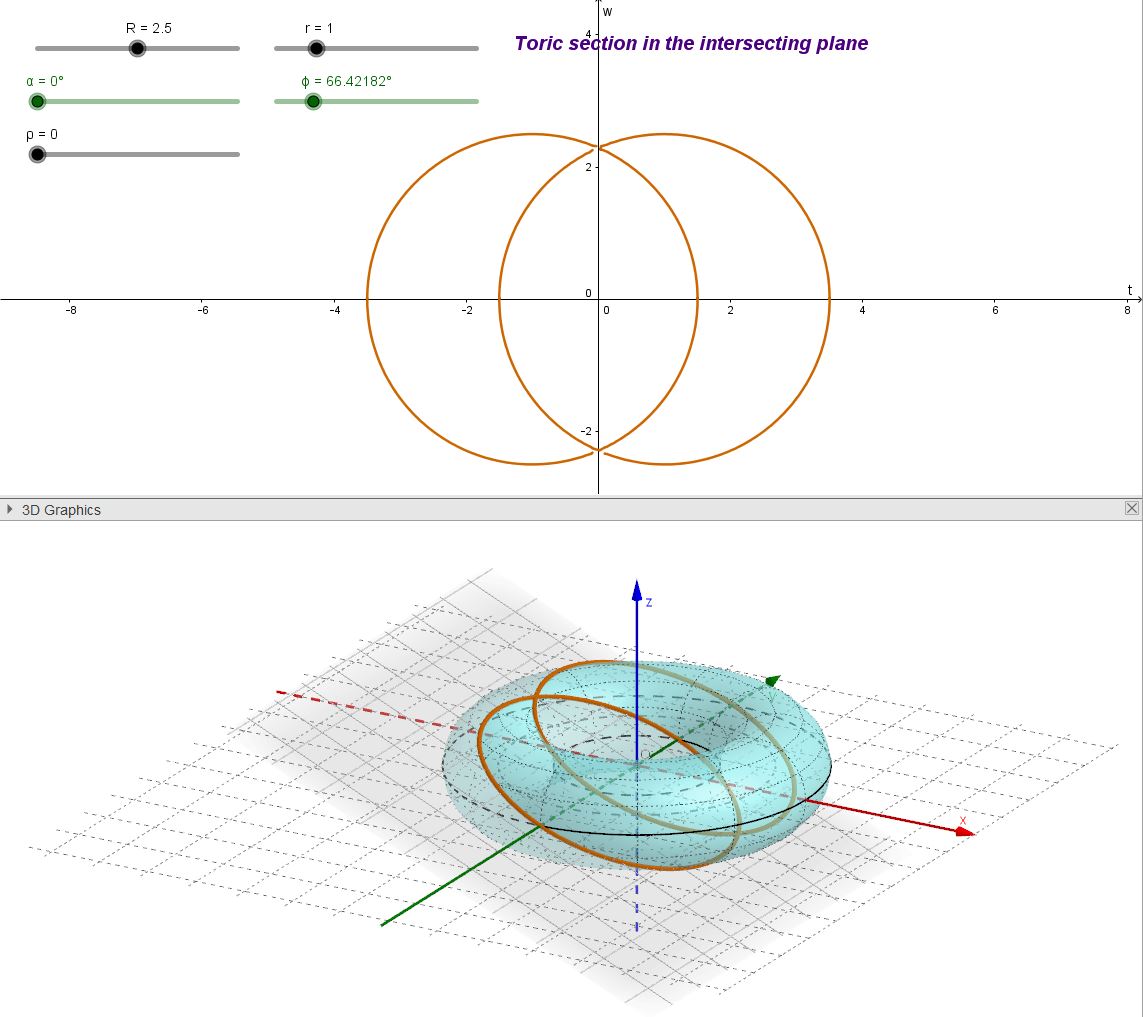}
	\end{minipage}
	\begin{minipage}{.49\textwidth}
	\centering
	\includegraphics[width=\linewidth]{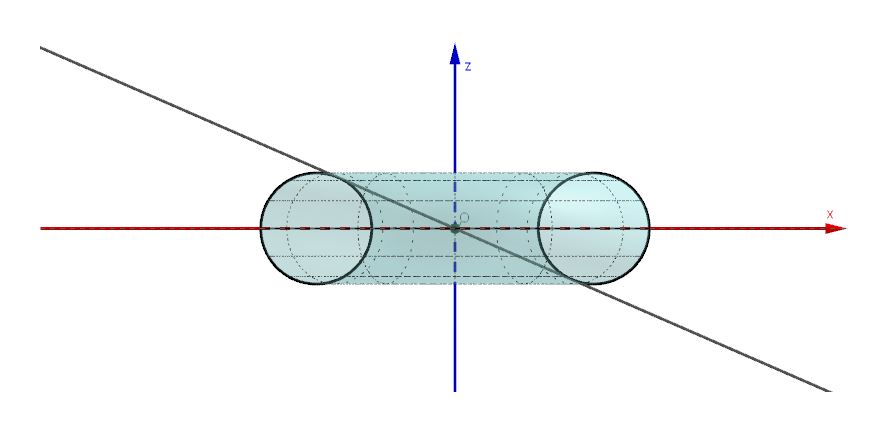}
	\end{minipage}
	\caption{The Villarceau's circles and the position of the intersecting plane to generate the Villarceau's circles}
	\label{fig03}
\end{figure}

\subsection{Central toric sections --- other sections}
It must be noted that if some central conic section generates circles, not all central toric sections are circles, as in the following examples.

\begin{figure}[H]
	\centering
	\begin{minipage}{.32\textwidth}
	\centering
	\includegraphics[width=\linewidth]{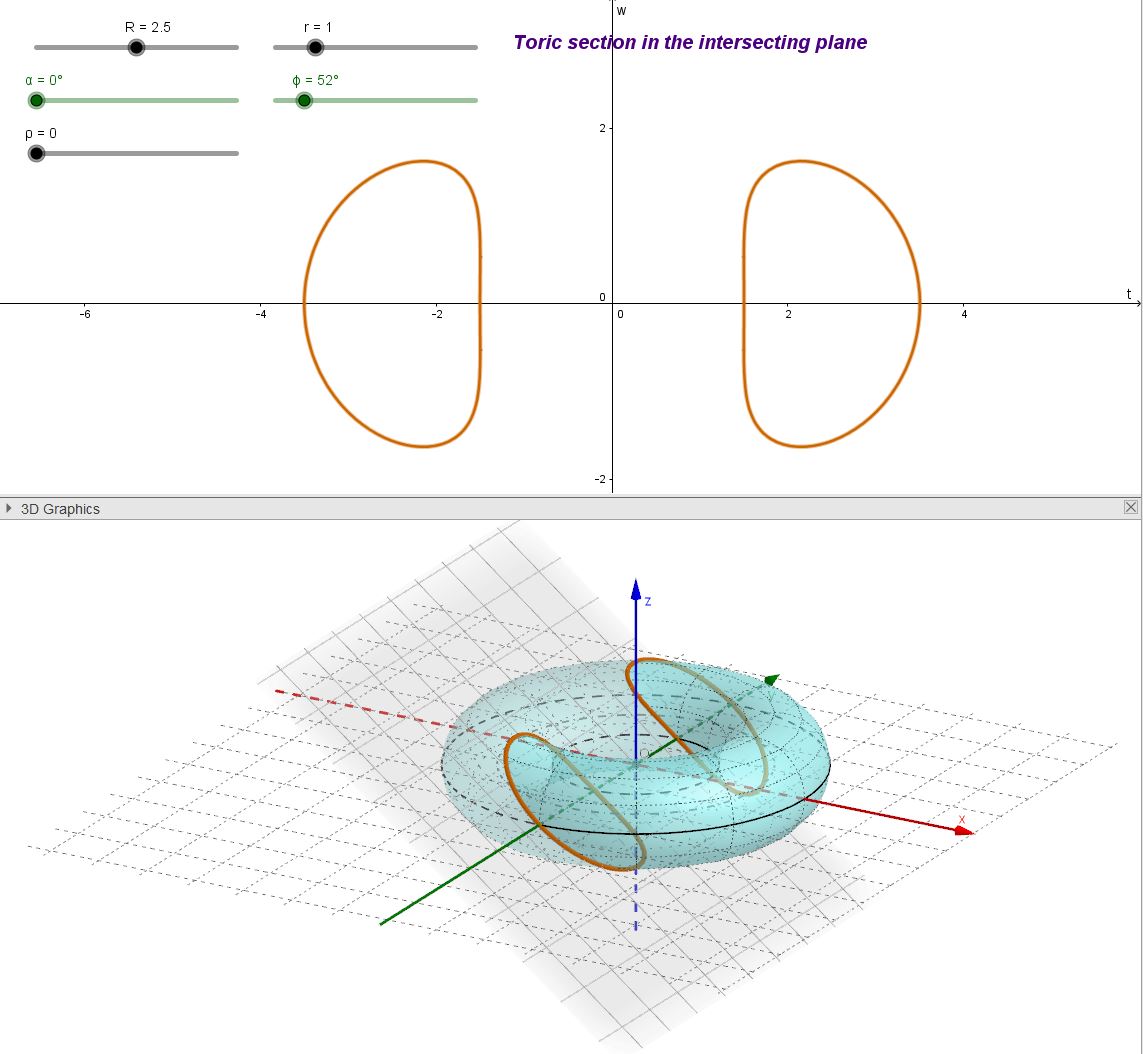}
	\end{minipage}
	\begin{minipage}{.32\textwidth}
		\centering
		\includegraphics[width=\linewidth]{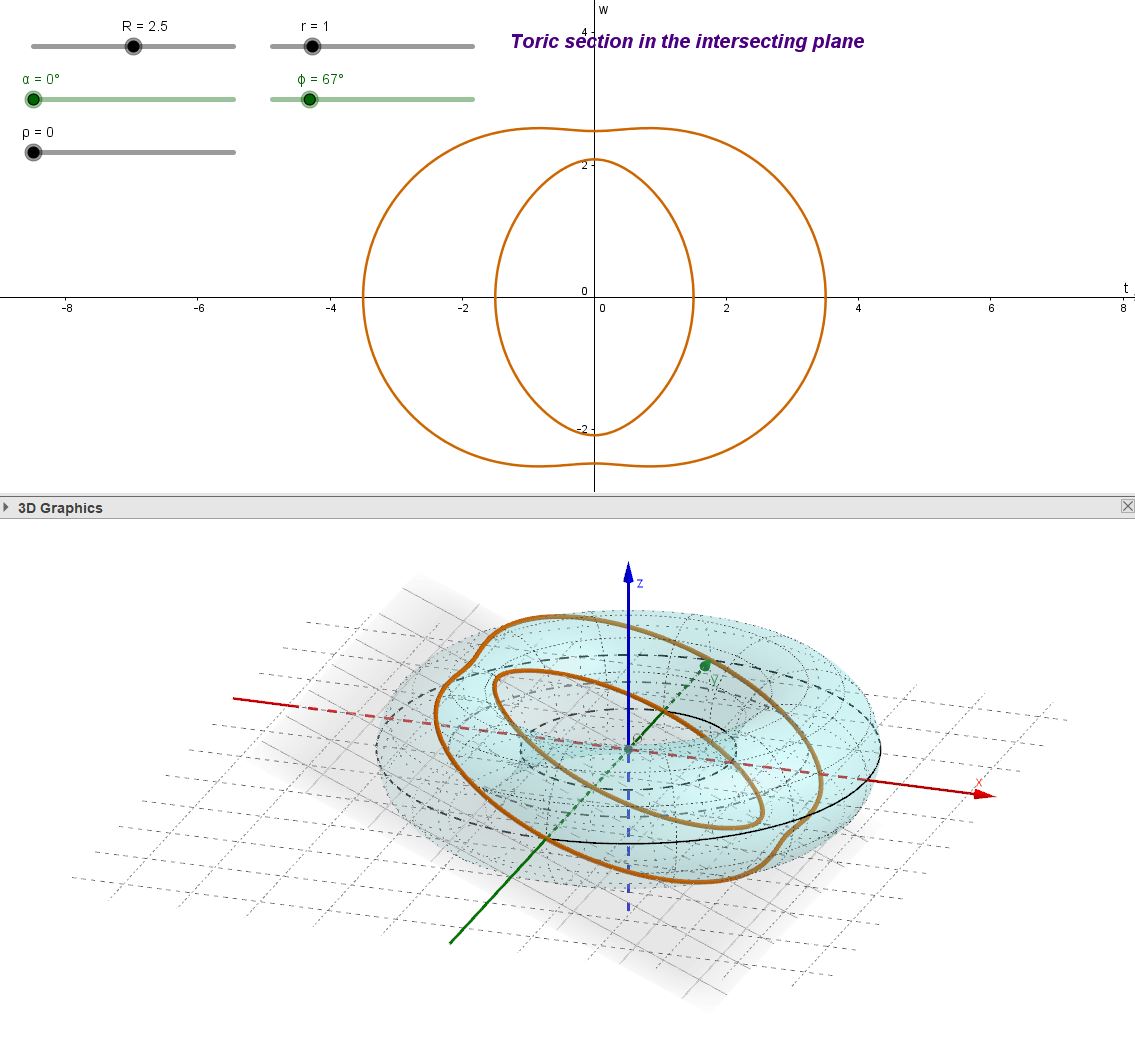}
	\end{minipage}
	\begin{minipage}{.32\textwidth}
		\centering
		\includegraphics[width=\linewidth]{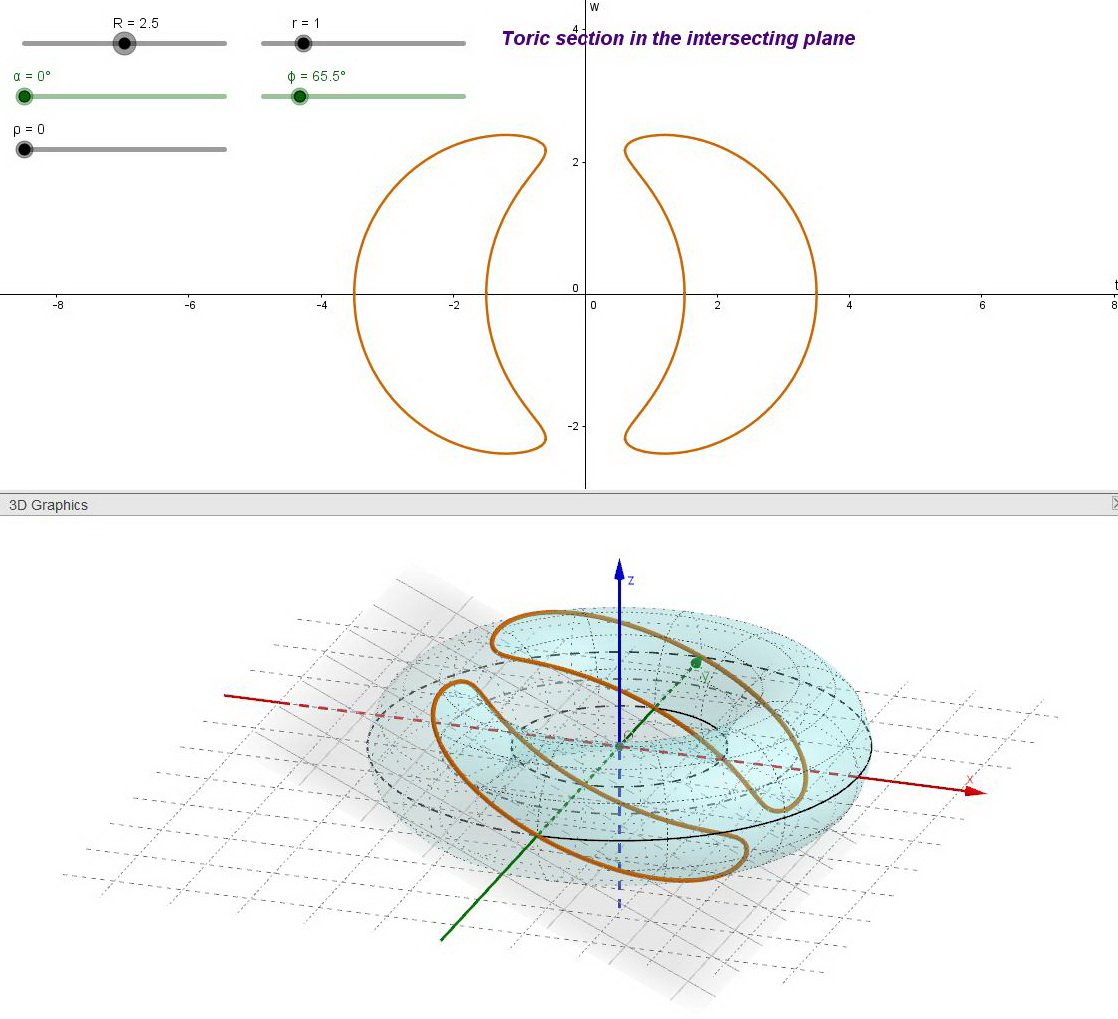}
	\end{minipage}
	\caption{Generic central toric sections}
	\label{fig04}
\end{figure}

\subsection{Spiric sections}
The toric section generated by an intersecting plane that is parallel to the torus axis (or perpendicular to the torus equatorial plane) is called a \textbf{\emph{spiric section}}. The name comes from the ancient Greek word \emph{$\sigma \pi \epsilon \iota \rho \alpha$} for torus (\cite{pac}).
Cassini's ovals and Bernoulli's Lemniscates belongs to this category.

\begin{figure}[H]
	\centering
	\begin{minipage}{.49\textwidth}
		\centering
		\includegraphics[width=\linewidth]{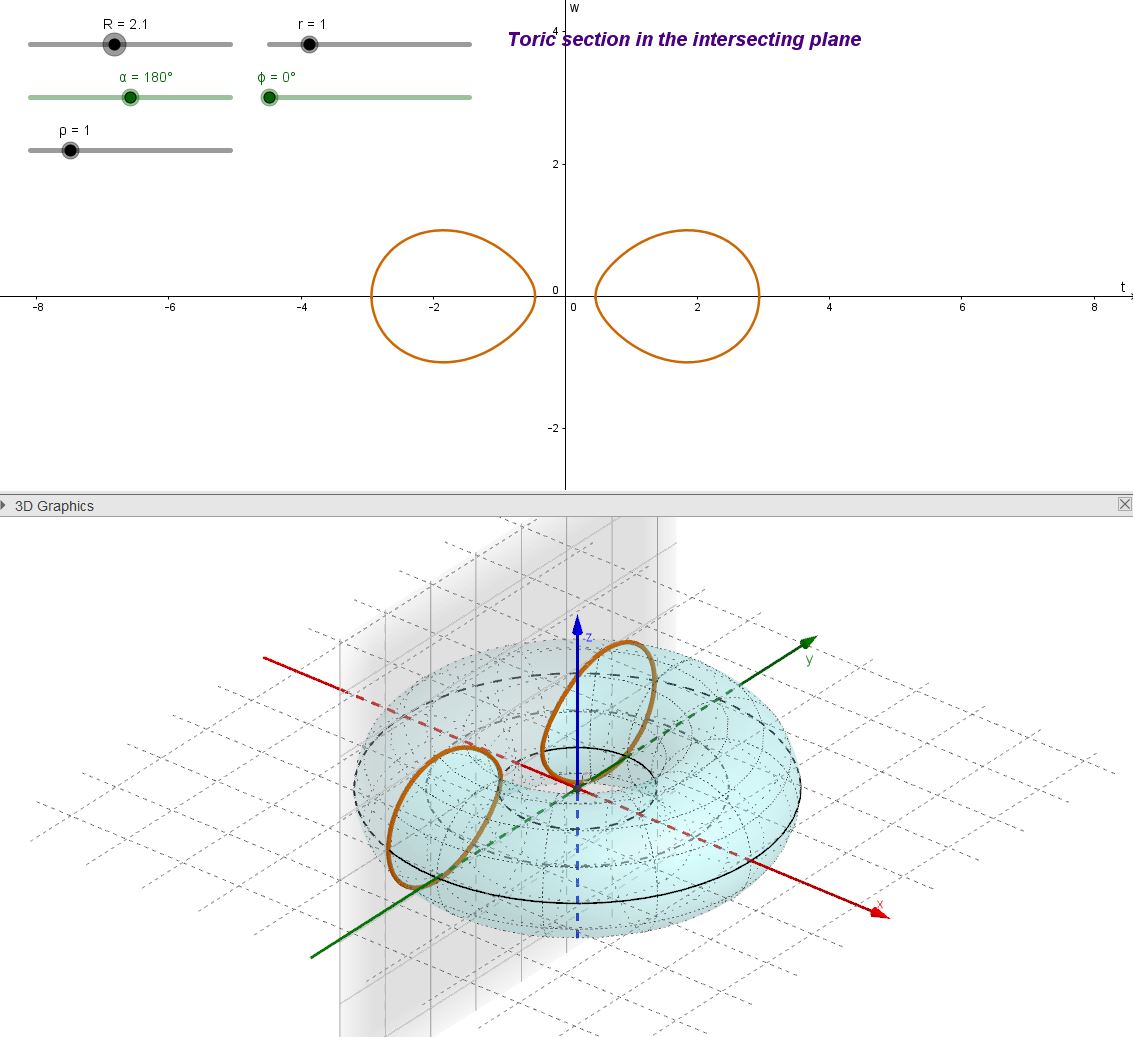}
	\end{minipage}
	\begin{minipage}{.49\textwidth}
		\centering
		\includegraphics[width=\linewidth]{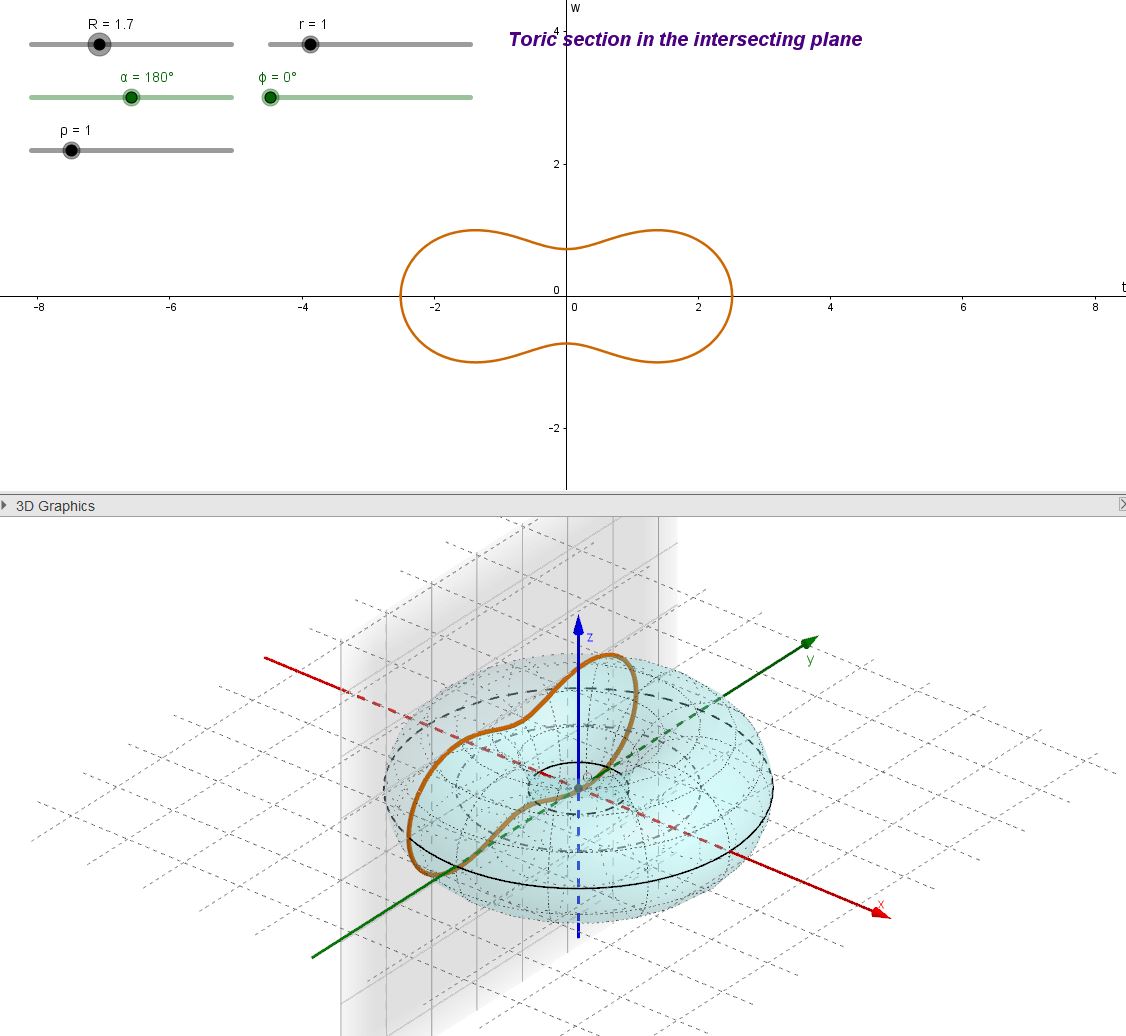}
	\end{minipage}
	\caption{Two different Cassini's ovals}
	\label{fig05}
\end{figure}

\textbf{Cassini's ovals} (Figure \ref{fig05}) are spiric sections in which the distance $\rho$ of the cutting plane to the torus axis equals the radius $r$ of the generating circle.
\par
An interesting alternative definition of a Cassini's oval is that of the set of points $ P $ such that the product of their distances to two fixed points $ F_1 $ and $ F_2 $ is constant:
$$\overline{P{{F}_{1}}}\cdot \overline{P{{F}_{2}}}={{b}^{2}}$$
This definition recalls that of the ellipse but here we must keep constant the \textit{product} of the distances instead of the \textit{sum}. 
If the focal distance is $\overline{{{F}_{1}}{{F}_{2}}}=2c$, the Cassini's ovals parameters, with reference to the torus intersected, have the values ${{b}^{2}}=2Rr$ and $c=R$.

\begin{figure}[H]
	\centering
	\includegraphics[width=3.7in]{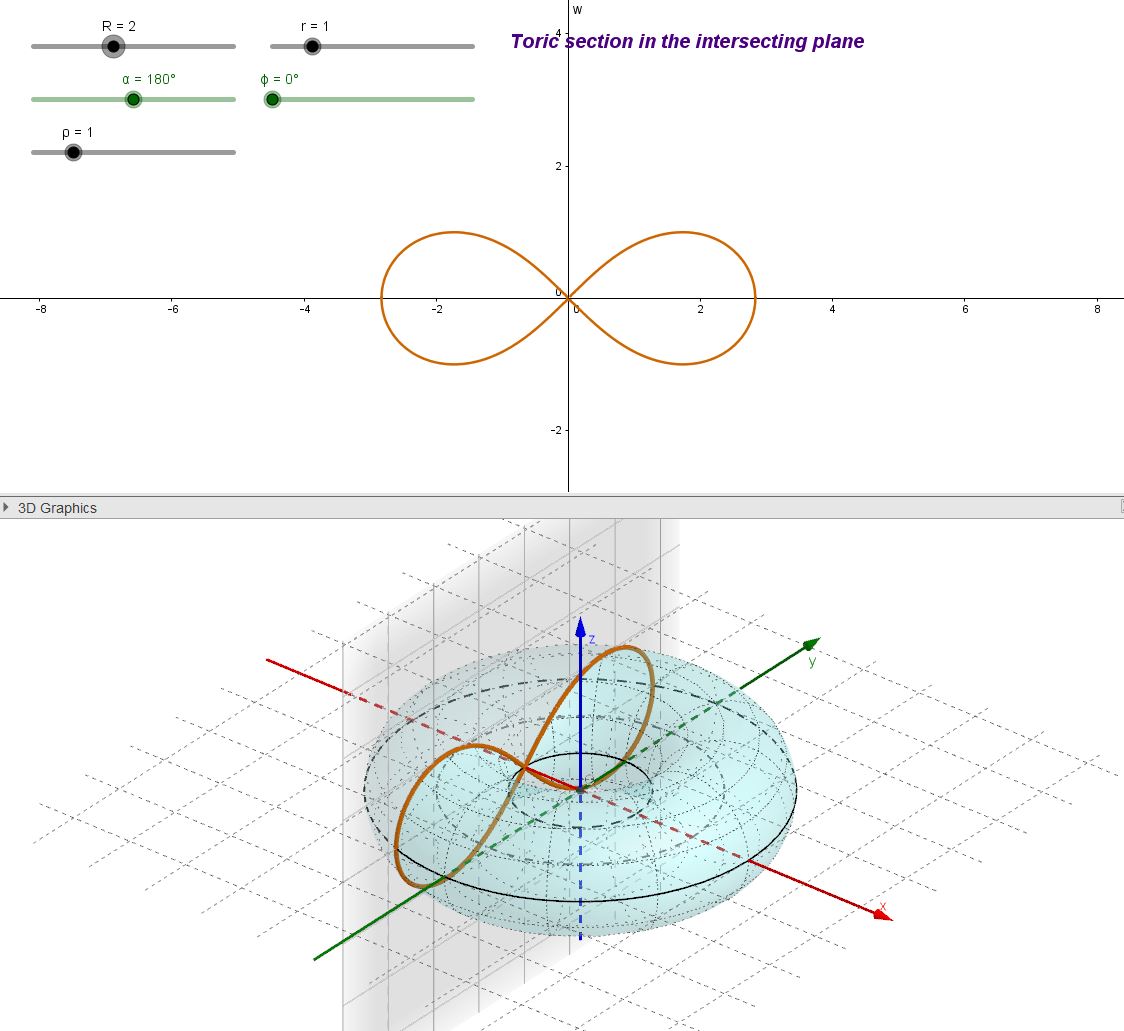}
	\caption{Bernoulli's lemniscate}
	\label{fig06}
\end{figure}

\textbf{Bernoulli's lemniscate} (Figure \ref{fig06}) is, in turn, a special case of a Cassini's oval. It is generated as a Cassini's oval with the further requirement that $R=2r$. In practice the intersecting plane must be tangent to the inner equator of the torus.

The lemniscate of Bernoulli can also be defined as the set of points $P$ such that the product of their distances to two fixed points $ F_1 $ and $ F_2 $ is constant and is also equal to the semi focal distance squared:
$$\overline{P{{F}_{1}}}\cdot \overline{P{{F}_{2}}}={{c}^{2}} {\quad\rm{where}\quad} \overline{{{F}_{1}}{{F}_{2}}}=2c$$ 

As a curiosity we can note that this lemniscate curve, representing an "\textit{eight}" rotated by $90^\circ$, has been used as the symbol for mathematical "infinity" ($\infty$) starting from the 17th century\footnote{It seems that the first appearance of the infinity symbol with its actual mathematical meaning was in the tratise "\textit{De sectionibus conicis}" by John Wallis in 1655.}

\begin{figure}[H]
	\centering
	\includegraphics[width=4.0in]{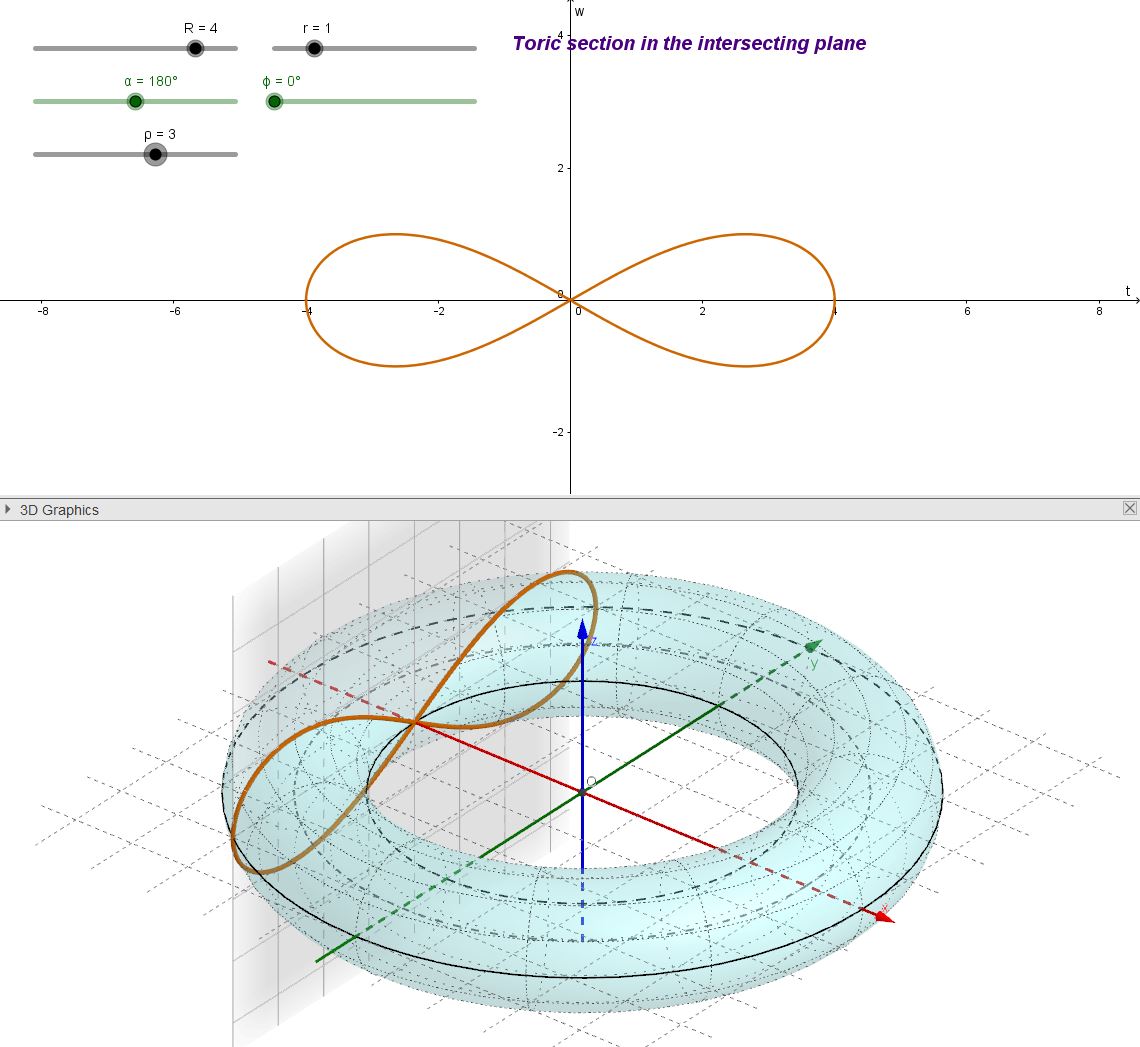}
	\caption{Hippopede of Proclus}
	\label{fig07}
\end{figure}

If we drop the requirement $\rho =r$ that characterizes Cassini's ovals but keep the condition that the intersecting plane is tangent to the interior circle (so that it must be $\rho =R-r$) we have the family of curves called Hippopedes of Proclus (Figure \ref{fig07}).

\clearpage
\subsection{Generic toric sections}

\begin{figure}[H]
	\centering
	\includegraphics[width=0.7 \linewidth]{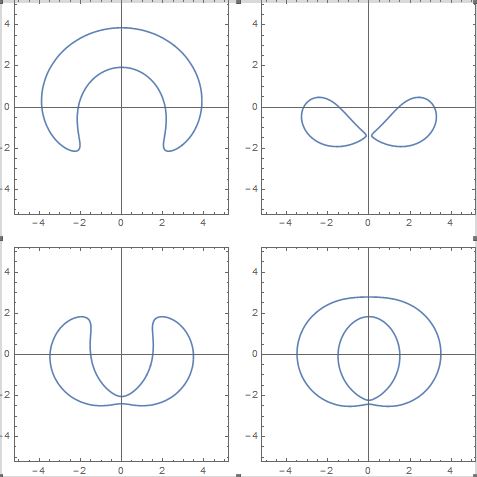}
	\caption{More generic toric sections}
	\label{fig08}
\end{figure}

There are many other toric sections that are neither central sections nor spiric sections. This means that the cutting plane is generically slanted with respect to the torus equatorial plane and that it has a non-zero distance from the center of the torus  (Figure \ref{fig08}). 

\section{Toric sections and conic sections}
The torus surface \ref{fig10} is generated by rotating a generating circle around an axis coplanar with the circle itself. In other words it is produced by moving around a first circle of center $O$ the center $O'$ of a second circle, providing that the second circle plane stays perpendicular to the tangent line to the first circle in $O'$.

On the other hand a conical surface (\textit{right circular conical surface}) \ref{fig09} can be generated by rotating a straight line around another straight line, providing that the first intersect the second, keeping fixed the angle between the two lines. 

The main difference between the two surfaces can then summarized as "\textit{circle rotated around a line}" vs. "\textit{line rotated around line}".

Since in a Cartesian coordinate system a circle is described by a second order equation and a line by a first order equation and that the rotation involves following a circular path, it can be anticipated that, in a Cartesian coordinate system, the equation of a toric section can be a quartic whilst, for a conic section, a quadratic equation suffices.

\section{The conical surface equation and the torus equation}
A \textbf{right circular conical surface} of aperture $2\theta$, whose axis is the $z$ coordinate axis, and whose apex is the origin $O$, is described by the parametric equations:
$${C_\theta }\left( {s,\psi } \right) = \left\{ {\begin{array}{*{20}{l}}
	{x = s\sin \theta \cos \psi }\\
	{y = s\sin \theta \sin \psi }\\
	{z = s\cos \theta }
	\end{array}} \right.$$
where, if $P$ is a point on the surface, the parameter $s$ is the coordinate of $P$ along the oriented line $OP$ and the parameter $\psi$ is the angle of rotation of $P$ around the $z$ axis.

\begin{figure}[H]
	\centering
	\includegraphics[width=0.7 \linewidth]{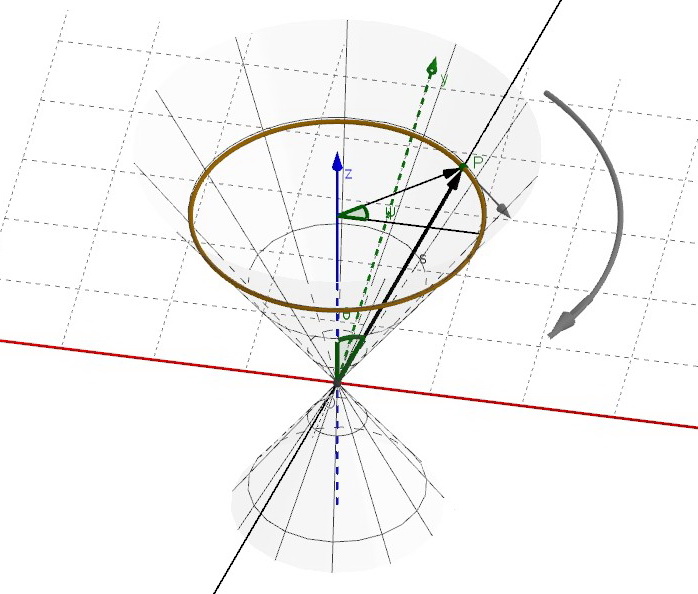}
	\caption{The conical surface}
	\label{fig09}
\end{figure}

Figure \ref{fig09} explains the meaning of the parametric equations.

In Cartesian coordinates the conical surface has implicit equation:
$$\left( {{x^2} + {y^2}} \right){\cos ^2}\theta  - {z^2}{\sin ^2}\theta  = 0$$

\noindent A \textbf{torus} centered in the origin with the $z$ axis as axis of revolution and with $R$ and $r$ as the major and minor radii, is described by the parametric equations:
$${T_{R,r}}\left( {u,v} \right)\left\{ {\begin{array}{*{20}{l}}
	{x = \left( {R + r\cos u} \right)\cos v}\\
	{y = \left( {R + r\cos u} \right)\sin v}\\
	{z = r\sin v}
	\end{array}} \right.$$
where, if $P$ is a point on the surface, the parameter $u$ is the angular position of $P$ on the revolving circle and $v$ is the angular position of the plane of the revolving circle with respect to the $z$ axis.

\begin{figure}[H]
	\centering
	\includegraphics[width=0.8 \linewidth]{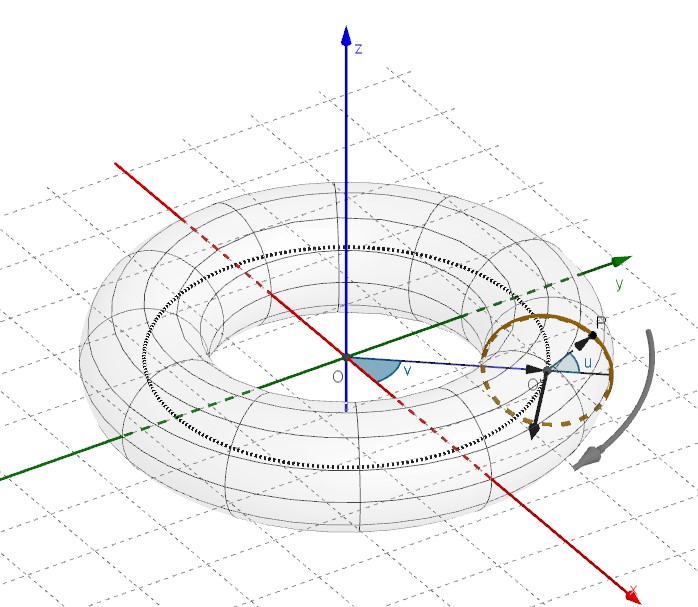}
	\caption{The torus}
	\label{fig10}
\end{figure}

Figure \ref{fig10} explains the meaning of the parametric equations.

\par \noindent The torus Cartesian (implicit) equation:
$${\left( {\sqrt {{x^2} + {y^2}} - R} \right)^2} + {z^2} = {r^2}$$
or, eliminating the square root to obtain a polynomial equation  in $x$, $y$, $z$\footnote{This equation belongs to the family of bicircular quartics \cite{bicqua}, that have the general form $${{\left( {{x}^{2}}+{{y}^{2}} \right)}^{2}}+a{{x}^{2}}+b{{y}^{2}}+cx+dy+e=0$$}
$${\left( {{x^2} + {y^2} + {z^2} + {R^2} - {r^2}} \right)^2} - 4{R^2}\left( {{x^2} + {y^2}} \right) = 0$$

\section{Analogies and differences between the conical surface and the torus (and their equations)}

\subsection{Analogies}
Both surfaces are bi-dimensional geometric objects immersed in a three dimensional space and so, like for any other surface, their parametric equations (if they can be made explicit), are always controlled by two variable parameters that set the coordinates of every point of the surface.

Both surfaces can be generated as surface of revolution, rotating a simpler figure around an axis.

Both surfaces have a very high degree of symmetry, as they are symmetric not only with respect to the origin, but also with respect of each of the three Cartesian coordinate axes (and this implies that they are also symmetric with respect to the three Cartesian planes). They are also symmetric with respect to all planes that are perpendicular to the horizontal plane and are running through the origin.

This high level of symmetry follows form the simple fact that the Cartesian equations of both surfaces only contains the squares of the three Cartesian coordinates, so that changing, say, $x$ with ($-x$) has no effect on the equations.

\subsection{Differences}
The conical surface is unbounded and is characterized by one single parameter determining its shape (the angle $\theta $).

The torus has a finite extension and is characterized by two parameters determining its shape (the major radius $R$ and the minor radius $r$) and size, but, actually its "scale free" shape is governed by the $R/r$ ratio.

\subsection{Topological differences (quick note for the more intrepid readers)}
Topologically a torus is the product space (or Cartesian product) of two circles, that can be expressed as
$${{\bf{T}}^{\bf{2}}} = {S^1} \times {S^1}$$

The topological definition of a conical surface is a little more complicated.
The cone CX of a topological space X is the quotient space:
$$CX = (X \times I)/(X \times \{ 0\} )$$
of the product of X with the unit interval $I = [0, 1]$.
Intuitively we make $X$ into a cylinder and collapse (identify or glue) one end of the cylinder to a point.

More informally a torus can be thought of as a cylindrical surface with the top and bottom circles glued together, whilst a conical surface can be seen as a cylindrical surface in which one of its circular section is collapsed to a point.

\subsection{Analogies and differences in the equations}
In the parametric equations of both surfaces there are trigonometric functions (to represent the rotations), but in the equation of the conical surface one of the parameters ($s$) acts linearly.

The Cartesian equations of both are given in implicit form and the Cartesian coordinates always appears with an even grade. The polynomial form of the conical surface is of second order, while the polynomial form of the torus is fourth order.

\section{Building the equation of a toric section}
Before attacking the problem of finding the equation of the curve of intersection between a torus and a plane it's necessary to examine how a plane can be described by an equation and which form (Cartesian or parametric) is more convenient for the purpose.

\subsection{The plane}
The plane can be identified by the vector orthogonal to it that starts at the origin and ends in the projection point of the origin on the plane. So let's use the parameters $\alpha$, $\phi$ and $\rho$ that set the position of the normal vector $\vec u$ through the origin (fig. 11).

$\alpha$ is the direction in the horizontal $xy$ plane (azimuthal angle) $\phi$ is the direction with respect of the horizontal plane (elevation angle) and $\rho$ is the modulus (length) of the vector.
The components of the vector are then $\vec u = \left( {\rho \cos \alpha \cos \phi ,\rho \sin \alpha \cos \phi ,\rho \sin \phi } \right)$.

These are also the coordinates of the point $Q$, projection of the origin $O$ on the plane.
The Cartesian equation of the plane is then:
$$p:\,\,\,\left( {x - {x_Q}} \right)\cos \alpha \cos \phi + \left( {y - {y_Q}} \right)\sin \alpha \cos \phi + \left( {z - {z_Q}} \right)\rho \sin \phi = 0$$
and its alternative parametric equations are:
$$p:{\kern 1pt} {\kern 1pt} {\kern 1pt} \left\{ {\begin{array}{*{20}{l}}
	{x = {x_Q} + t\sin \alpha  - w\cos \alpha \sin \phi }\\
	{y = {y_Q} - t\cos \alpha  - w\sin \alpha \sin \phi }\\
	{z = {z_Q} + w\cos \phi }
	\end{array}} \right.$$

\begin{figure}[H]
	\centering
	\includegraphics[width=0.7 \linewidth]{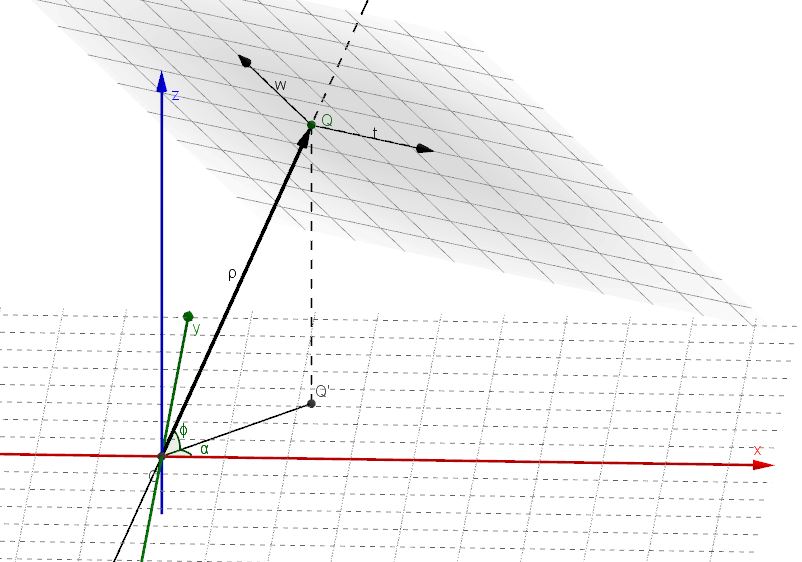}
	\caption{The plane in tridimensional space}
	\label{fig11}
\end{figure}

With above parametric equations the parameters $t$ and $w$  can be interpreted as the \textit{embedded} orthogonal Cartesian axes in the plane starting from the point $Q$ (plane origin), where $t$ is the horizontal axis (parallel to the $xy$ plane) and $w$ is the vertical axis (perpendicular to $t$) in this plane.

\subsection{The intersection curve between the plane and the torus}
We can now use the following equations of the torus and the plane to get the intersection curve:
$$T\left( {x,y,z} \right) = {\left( {\sqrt {{x^2} + {y^2}}  - R} \right)^2} + {z^2} - {r^2} = 0$$
$$p\left( t,w \right):\left\{ \begin{array}{*{35}{l}}
x={{x}_{Q}}+t\sin \alpha-w\cos \alpha \sin \phi   \\
y={{y}_{Q}}-t\cos \alpha-w\sin \alpha \sin \phi   \\
z={{z}_{Q}}+w\cos \phi   \\
\end{array} \right.	\rm{\quad where \quad}\left\{ {\begin{array}{*{20}{l}}
	{{x_Q} = \rho \cos \alpha \cos \phi }\\
	{{y_Q} = \rho \sin \alpha \cos \phi }\\
	{{z_Q} = \rho \sin \phi }
	\end{array}} \right.$$

We'll plug in the $x$, $y$, $z$ parametric equations of the plane in the implicit equation of the torus. In so doing we'll then get an implicit equation in the new variables $t$ and $w$ and since these variables are just the “$xy$” coordinates inside the intersecting plane, this implicit equation will be nothing but the implicit equation of the toric section as seen in the intersecting plane.

Let's start
\footnotesize $${\left( {\sqrt {{{\left( {{x_Q} + t\sin \alpha {\rm{\;}} - w\cos \alpha \sin \phi } \right)}^2} + {{\left( {{y_Q} - t\cos \alpha {\rm{\;}} - w\sin \alpha \sin \phi } \right)}^2}} {\rm{\;}} - R} \right)^2} + {\left( {\rho \sin \phi  + w\cos \phi } \right)^2} - {r^2} = 0$$ \normalsize
Expanding the squares, the expression inside the root simplifies to
${t^2} + {\left( {\rho \cos \phi  - w\sin \phi } \right)^2}$
so our equation will be:
$${\left( {\sqrt {{t^2} + {{\left( {\rho \cos \phi  - w\sin \phi } \right)}^2}}  - R} \right)^2} + {\left( {{z_Q} + w\cos \phi } \right)^2} - {r^2} = 0$$
Expanding the squares once again we finally get
$${t^2} + {w^2} + {\rho ^2} + {R^2} - {r^2} = 2R\sqrt {{t^2} + {{\left( {\rho \cos \phi  - w\sin \phi } \right)}^2}}$$
that is, renaming the ($t$, $w$) orthogonal coordinates in the plane with the usual letters ($x$, $y$) the equation will be
$${x^2} + {y^2} + {\rho ^2} + {R^2} - {r^2} = 2R\sqrt {{x^2} + {{\left( {\rho \cos \phi  - y\sin \phi } \right)}^2}}$$
or, in polynomial form,
$${\left( {{x^2} + {y^2} + {\rho ^2} + {R^2} - {r^2}} \right)^2} - 4{R^2}\left( {{x^2} + {{\left( {\rho \cos \phi  - y\sin \phi } \right)}^2}} \right) = 0$$
hence
\footnotesize $${{\left( {{x^2} + {y^2}} \right)^2} + {\left( {{\rho ^2} + {R^2} - {r^2}} \right)^2} + 2\left( {{x^2} + {y^2}} \right)\left( {{\rho ^2} + {R^2} - {r^2}} \right) - 4{R^2}\left( {{x^2} + {\rho ^2}{{\cos }^2}\phi  + {y^2}{{\sin }^2}\phi  - 2y\rho \cos \phi \sin \phi } \right) = 0}$$ \normalsize
Above equation is in the form:
$${\left( {{x^2} + {y^2}} \right)^2} + a{x^2} + b{y^2} + cy + d = 0$$
It can be noted that this equation is symmetric with respect to the $x$ axis (but not with respect to the $y$ axis).

\subsection{The equation of the toric section in space}
Having found an implicit equation giving the toric section represented in the intersecting plane, the last problem is to translate this equation in space and get an equation presenting the same curve in the tridimensional space.
Geogebra can draw even complicated curves in the 3D graphics view if they are expressed in a parametric form, with a single parameter, of the type:
$$t\left( w \right):\left\{ \begin{array}{*{35}{l}}
x={{f}_{x}}\left( w \right)  \\
y={{f}_{y}}\left( w \right)  \\
z={{f}_{z}}\left( w \right)  \\
\end{array} \right.$$
Anyway, the implicit equation previously found contains two parameters: $t$ and $w$. Furthermore we know that the curve lies on the intersecting plane.
What can be done is find the relation between $t$ and $w$ given by the implicit equation, expressing one parameter as a function of the other, and use the plane equation with the second parameter.
The examination of the implicit equation
$${\left( {\sqrt {{t^2} + {{\left( {\rho \cos \phi  - w\sin \phi } \right)}^2}}  - R} \right)^2} + {\left( {{z_Q} + w\cos \phi } \right)^2} - {r^2} = 0$$
suggests to express $t$ as a function of $w$: $t = t\left( w \right)$
With some algebra we get to
$${t = t\left( w \right) =  \pm \sqrt { - {{\left( {\rho \cos \phi  - w\sin \phi } \right)}^2} + {{\left( {R \pm \sqrt {{r^2} - {{\left( {w\cos \phi  + \rho \sin \phi } \right)}^2}} } \right)}^2}} }$$
Not much meaningful and simple, we admit...
Our curve in space will then be
$$\tau \left( w \right):\left\{ \begin{array}{*{35}{l}}
x={{x}_{Q}}+t\left( w \right)\sin \alpha-w\cos \alpha \sin \phi   \\
y={{y}_{Q}}-t\left( w \right)\cos \alpha-w\sin \alpha \sin \phi   \\
z={{z}_{Q}}+w\cos \phi   \\
\end{array} \right.$$
and Geogebra can handle it and produce the toric section in tridimensional space.

\section{A surprise! (the conical bridge)}
The implicit equation of the toric section in the intersecting plane
$${{t}^{2}}+{{w}^{2}}+{{\rho }^{2}}+{{R}^{2}}-{{r}^{2}}=2R\sqrt{{{t}^{2}}+{{\left( \rho \cos \phi -w\sin \phi  \right)}^{2}}}$$
can be written as follows by renaming the variables ($t$, $w$) with the usual symbols ($x$, $y$)
$${{x}^{2}}+{{y}^{2}}+{{\rho }^{2}}+{{R}^{2}}-{{r}^{2}}=2R\sqrt{{{x}^{2}}+{{\left( \rho \cos \phi -y\sin \phi  \right)}^{2}}}$$
If we try to simplify it a little with some change of variables to get rid of the square root, we could have the idea to make the transformation
$$\left\{ \begin{array}{*{35}{l}}
{{x}^{2}}+{{\left( \rho \cos \phi -y\sin \phi  \right)}^{2}}={{k}^{2}}{{z}^{2}}  \\
y=y  \\
\end{array} \right.$$
that is, we set ${x^2}$ as
${{x}^{2}}={{k}^{2}}{{z}^{2}}-{{\left( \rho \cos \phi -y\sin \phi  \right)}^{2}}$

So we apply the geometric transformation that maps $\left( {{x}^{2}},y \right)$ to $\left( {{k}^{2}}{{z}^{2}}-{{\left( \rho \cos \phi -y\sin \phi  \right)}^{2}},y \right)$

With this mapping the toric section equation becomes
\footnotesize
$${{{k}^{2}}{{z}^{2}}-{{\left( \rho \cos \phi -y\sin \phi  \right)}^{2}}+{{y}^{2}}+{{\rho }^{2}}+{{R}^{2}}-{{r}^{2}}=2R\sqrt{{{k}^{2}}{{z}^{2}}-{{\left( \rho \cos \phi -y\sin \phi  \right)}^{2}}+{{\left( \rho \cos \phi -y\sin \phi  \right)}^{2}}}}$$ 
\normalsize

that is (in the $\left( z,y \right)$ space)
$${{k}^{2}}{{z}^{2}}-{{\left( \rho \cos \phi -y\sin \phi  \right)}^{2}}+{{y}^{2}}+{{\rho }^{2}}+{{R}^{2}}-{{r}^{2}}=2R\sqrt{{{k}^{2}}{{z}^{2}}}$$
$${{k}^{2}}{{z}^{2}}-{{\rho }^{2}}{{\cos }^{2}}\phi -{{y}^{2}}{{\sin }^{2}}\phi +2y\rho \cos \phi \sin \phi +{{y}^{2}}+{{\rho }^{2}}+{{R}^{2}}-{{r}^{2}}=2R\left( \pm kz \right)$$
$${{k^2}{z^2} + {y^2}\left( {1 - {{\sin }^2}\phi } \right) - 2R\left( { \pm z} \right) + 2y\rho \cos \phi \sin \phi  + {\rho ^2}\left( {1 - {{\cos }^2}\phi } \right) + {R^2} - {r^2} = 0}$$
$${{k}^{2}}{{z}^{2}}\pm 2Rkz+{{R}^{2}}+{{y}^{2}}{{\cos }^{2}}\phi +2y\rho \cos \phi \sin \phi +{{\rho }^{2}}{{\sin }^{2}}\phi ={{r}^{2}}$$
$${{k}^{2}}{{\left( z\pm \frac{R}{k} \right)}^{2}}+{{\cos }^{2}}\phi {{\left( y+\rho \tan \phi  \right)}^{2}}={{r}^{2}}$$
Luckily we saw that many simplifications could be made. Good.
We have obtained the equation of two ellipses in the $yz$ plane. But we can do even better and make them two circles bysetting  $k=\cos \phi$, so it is
$${{\left( z\pm \frac{R}{\cos \phi } \right)}^{2}}+{{\left( y+\rho \tan \phi  \right)}^{2}}=\frac{{{r}^{2}}}{{{\cos }^{2}}\phi }$$
Above is the equation of a pair of circles in the $yz$ plane.\\
The centers are ${C_{1,2}}\left( { \pm \frac{R}{{\cos \phi }}; - \rho \tan \phi {\rm{\;}}} \right)$. The common radius is ${\frac{r}{\cos \phi }}$.

The meaning of all this is that if we start with the following two circles in a yz plane
$${{\left( z\pm \frac{R}{\cos \phi } \right)}^{2}}+{{\left( y+\rho \tan \phi  \right)}^{2}}=\frac{{{r}^{2}}}{{{\cos }^{2}}\phi }$$
and apply the transformation $\left( z,y \right)\mapsto \left( x,y \right)$ defined by
$$\chi :\left\{ \begin{array}{*{35}{l}}
{{x}^{2}}={{z}^{2}}{{\cos }^{2}}\phi -{{\left( \rho \cos \phi -y\sin \phi  \right)}^{2}}  \\
y=y  \\
\end{array} \right.$$
we get a toric section!

Now what is the geometric meaning of $\chi$?
The $y$ variable is unchanged, but the $z$ variable is related to the variables ($x$, $y$) by the equation
$${{z}^{2}}={{\left( \frac{x}{\cos \phi } \right)}^{2}}+{{\left( \rho -y\tan \phi  \right)}^{2}}$$

Above is, in tridimensional space, the Cartesian equation of a conical surface with vertex ${V\left( 0;\frac{\rho }{\tan \phi };0 \right)}$ and a vertical axis.

This is not a circular conical surface and the sections of this surface with the horizontal planes $z=k$ are ellipses with semi axes of length $k\cos \phi$ and $\frac{k}{\tan \phi }$.

\begin{figure}[H]
	\centering
	\includegraphics[width=0.7 \linewidth]{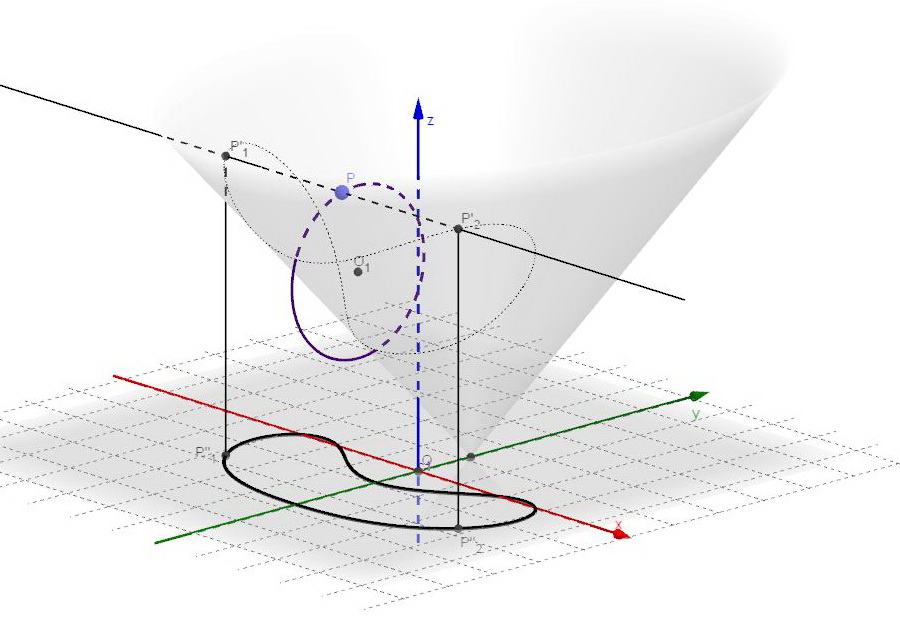}
	\caption{Construction of a toric section with a conical surface}
	\label{fig12}
\end{figure}

So, to summarize the previous passages, here's a curious unexpected way to trace a toric section.
\begin{enumerate}
	\item Start form one of the two circles in the $yz$ plane (one is enough) and choose a point on it $P\left( {{y}_{P}};{{z}_{P}} \right)$
	\item Find the intersection between the line through $P$ perpendicular to the $yz$ plane and the conical surface $C:{x^2} = {z^2}{\cos ^2}\phi  - {\left( {\rho \cos \phi  - y\sin \phi {\rm{\;}}} \right)^2}$
	In general, if the line actually intersect the cone, there will be two intersection points: $$P{{\rm{'}}_{{\rm{1}},2}} = \left( { \pm \sqrt {{z^2}{{\cos }^2}\phi  - {{\left( {\rho \cos \phi  - y\sin \phi {\rm{\;}}} \right)}^2}} ,{y_p},{z_p}} \right)$$
	\item Project this two points on the xy plane. They will be: 
	$$P{\rm{'}}{{\rm{'}}_{1,2}} = \left( { \pm \sqrt {{z^2}{{\cos }^2}\phi  - {{\left( {\rho \cos \phi  - y\sin \phi {\rm{\;}}} \right)}^2}} ,{y_p},0} \right)$$These two points are points of the toric section.
	\item By moving the point $P$ on the circle, the points $P'$ will trace a curve on the conical surface and the projection of this curve on the $xy$ plane is the toric section corresponding to the set of parameters $\left( R,r,\rho ,\phi  \right)$ used.
\end{enumerate}
Figure \ref{fig12} shows the construction method.

Actually we have done nothing but intersecting the conical surface with a cylinder and have verified that the projection of the curve of intersection on the $xy$ plane is a toric section.

\begin{figure}[H]
	\centering
	\includegraphics[width=0.7 \linewidth]{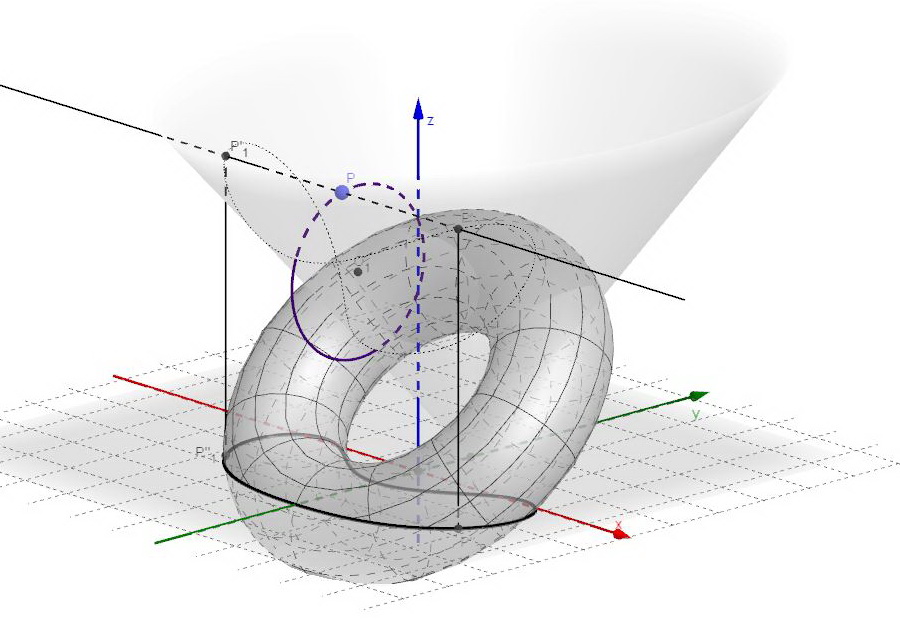}
	\caption{Check of the construction with the \textit{tilted} torus}
	\label{fig13}
\end{figure}

With this method we trace the toric section on the intersecting plane that, in this case, is the $xy$ plane.
So the torus intersected must have to be tilted by the angle $\frac{\pi }{2} - \phi$ to reproduce, in the new system of reference, the initial figure in which a horizontal torus was intersected by an oblique plane.
Figure \ref{fig13} shows the tilted torus intersecting the $xy$ plane. The previously obtained curve is actually a toric section.

\begin{figure}[H]
	\centering
	\includegraphics[width=0.65 \linewidth]{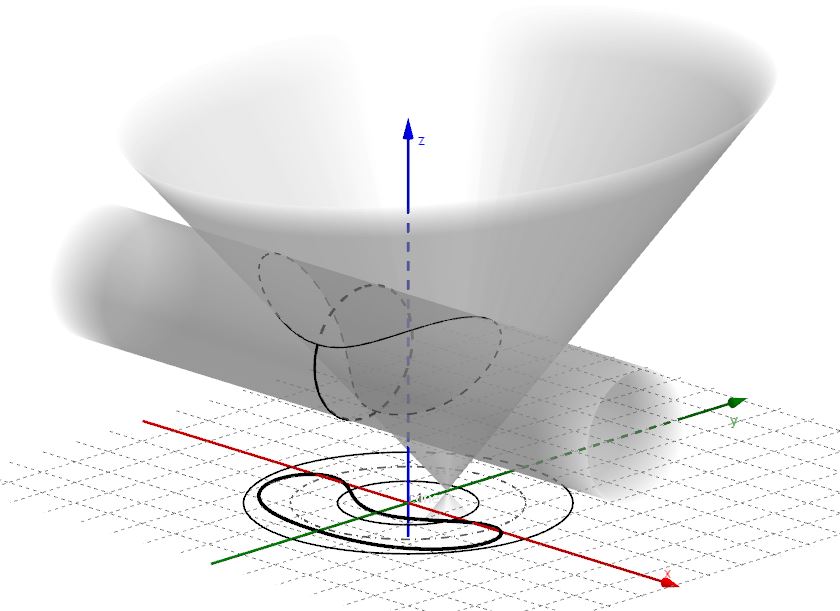}
	\caption{The toric section generated through a cone-cylinder intersection}
	\label{fig14}
\end{figure}

\begin{figure}[H]
	\centering
	\includegraphics[width=0.65 \linewidth]{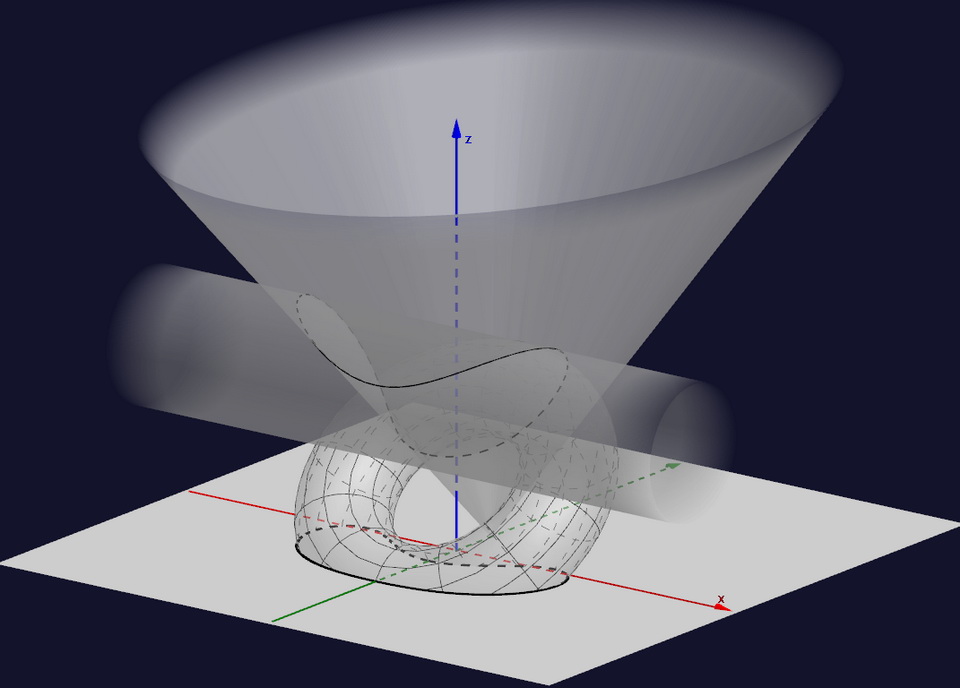}
	\caption{The torus-cylinder-cone interesting \textit{interplay}}
	\label{fig15}
\end{figure}

In conclusion, we have started with a comparison of toric and conic sections, derived the toric section equation (fourth grade), and, with some algebraic manipulation, found that the same toric section equation can also be seen has the projection on a plane of a cone-cylinder intersection (where both surfaces have second grade equations). The conical surface, then, is more closely related to the torus than expected as it enters, as a sort of bridge (Figures \ref{fig14} and \ref{fig15}), in the approximate abstract equivalence relation
$$torus\bigcap plane\approx cone\bigcap cylinder$$

\section{The toric section in Geogebra}
\noindent In the page \href{http://www.lucamoroni.it/simulations/intersection-torus-plane-simulation/}{http://www.lucamoroni.it/simulations/intersection-torus-plane-simulation/} there is an interactive simulation about generating toric sections together with some explanation of the interactive commands.

\noindent The Geogebra file can also be seen in Geogebra's material repository (\href{https://ggbm.at/MmTVuXYk}{https://ggbm.at/MmTVuXYk}) or, better, downloaded as a ".ggb"
(\href{https://www.geogebra.org/material/download/format/file/id/Vp95FtPQ}{file}) and run locally on the PC.

\noindent The free Geogebra {\em Classic desktop} program is available at the geogebra.org \href{https://www.geogebra.org/download}{download page}


\begin{thebibliography}{1}
\addcontentsline{toc}{section}{References}
	\bibitem{torus} \href{https://en.wikipedia.org/wiki/Torus}{https://en.wikipedia.org/wiki/Torus}
	
	\bibitem{toric section}  \href{https://en.wikipedia.org/wiki/Toric_section}{https://en.wikipedia.org/wiki/Toric\textunderscore section}
	
	\bibitem{conicalsurface} \href{https://en.wikipedia.org/wiki/Conical_surface}{https://en.wikipedia.org/wiki/Conical\textunderscore surface}
	
	\bibitem{conetopology} \href{https://en.wikipedia.org/wiki/Cone_(topology)}{https://en.wikipedia.org/wiki/Cone\textunderscore (topology)}

	\bibitem{greatestlove} Sym, Antoni, {\em "Darboux's greatest love"}: Phys. A: Math. Theor. 42 (2009) 404001.

	\bibitem{geomat} Geogebra material on toric sections: \href{https://www.geogebra.org/m/MmTVuXYk)}{\emph{Torus-plane intersection}}

	\bibitem{bicqua} \href{https://opus4.kobv.de/opus4-fau/files/2301/ThomasWernerDissertation.pdf}{Thomas Werner (2011)  {\em Rational families of circles and bicircular quartics} - pp. 118-125}

	\bibitem{pac} Brieskorn-Kn{\"o}rrer (1986), {\em Plane Algebraic Curves}, Birkhäuser - pp. 16-18
	

\end{thebibliography}
\end{document}